\documentclass{elsart3-1}

 \usepackage{epsfig}
 
\usepackage{amssymb,amsmath}

\usepackage[english,francais]{babel} 
 \usepackage{color}
\usepackage{caption,sidecap}

\usepackage{algorithm}
\usepackage{algorithmic}

\def\calN{{\mathcal N}}
\def\calD{{\mathcal D}}

\def\boldf{{\mathbf f}}

\definecolor{blck}{rgb}{0,0,0}
\definecolor{darkred}{rgb}{.6,.1,0}

\newcommand\nl[1]{\textcolor{blck}{#1}}
\newcommand\rev[1]{\textcolor{blck}{#1}}

\newtheorem{e-proposition}[theorem]{Proposition}

\newtheorem{e-definition}[theorem]{Definition\rm}

\renewcommand{\theequation}{\arabic{equation}}
\newcommand{\m}[1]{\mathbf{#1}}
\newcommand{\mD}{\m{D}}

\newcommand{\mL}{\m{L}}
\newcommand{\bu}{\mathbf{u}}
\newcommand{\vc}{\mathbf{c}}
\def\og{\leavevmode\raise.3ex\hbox{$\scriptscriptstyle\langle\!\langle$~}}
\def\fg{\leavevmode\raise.3ex\hbox{~$\!\scriptscriptstyle\,\rangle\!\rangle$}}

\def\calN{{\mathcal N}}
\def\R{\mathbb{R}}

\date{}

\begin{document}
\centerline{}
\begin{frontmatter}


\selectlanguage{english}

\vspace{-1.8in}
\title{Parametric Analytical Preconditioning and its Applications to the Reduced Collocation Methods}


\selectlanguage{english}

\author{Yanlai Chen$^{1}$},
\ead{yanlai.chen@umassd.edu}
\author{Sigal Gottlieb$^{1}$},
\ead{sgottlieb@umassd.edu}
\author{Yvon Maday$^{2,3}$},
\ead{maday@ann.jussieu.fr}

\address{1 - Department of
Mathematics, University of Massachusetts Dartmouth, 285 Old Westport Road, North Dartmouth, MA 02747, USA.
The research of the first author was partially supported by National Science Foundation grant DMS-1216928.
The research of the second author was partially supported by AFOSR grant 
FA-9550-12-1-0224.}
\address{2 - Sorbonne Universit\'es, UPMC Univ Paris 06, UMR 7598, Laboratoire Jacques-Louis Lions \& Institut Universitaire de France, F-75005, Paris, France}
\address{3 - Division of Applied Mathematics, Brown University, 182 George St, Providence, RI 02912, USA}

\medskip
\begin{center}
{
Presented by Olivier Pironneau}
\end{center}

\vspace{-0.5cm}


\begin{abstract}
\selectlanguage{english}
In this paper, we \nl{extend the recently developed reduced collocation method \cite{ChenGottlieb} 
to the nonlinear case, and} propose two analytical  preconditioning strategies. 
One is parameter independent and easy 
to implement, the other one  has the traditional affinity with respect to the parameters which allows for 
efficient implementation through an offline-online decomposition.
Overall, the preconditioning improves the quality of the error estimation uniformly on the parameter domain, and speeds up the convergence of the 
reduced solution to the truth approximation.

\vskip 0.5\baselineskip

\selectlanguage{francais}
\noindent{\bf R\'esum\'e} \vskip 0.5\baselineskip \noindent
On \' etend dans cette note la m\'ethode de collocation r\'eduite r\'ecemment introduite  dans  \cite{ChenGottlieb}  au cas non lin\'eaire et on propose deux strat\'egies de pr\'econditionnement  dont une est ind\'ependante des param\`etres et facile a mettre en oeuvre et l'autre poss\`ede la propri\'et\'e classique de d\'ecomposition affine qui permet une mise en oeuvre rapide en-ligne/hors-ligne. Ces strat\'egies am\'eliorent  la qualit\'e de l'approximation et la vitesse de convergence.

\end{abstract}
\end{frontmatter}

\selectlanguage{francais}
\vspace{-1.2cm}
\section*{Version fran\c{c}aise abr\'eg\'ee}

La m\'ethode de base r\'eduite classique (RBM)\cite{Barrett_Reddien,Noor_Peters,Peterson,Prudhomme_Rovas_Veroy_Maday_Patera_Turinici} pour l'approximation de la solution d'\'equations aux d\'eriv\'ees partielles (EDP) param\'etr\'ees du type $[\mathbb{L} (\mu)\, u_\mu] (x) = f(x; \mu), \,\, x \in \Omega \,\,\,\rev{\subset \mathbb{R}^n}$ repose sur la d\'efinition d'un espace de discr\'etisation ad'hoc, engendr\'e par des solutions particuli\`eres de l'EDP en certain param\`etres bien choisis. Ces solutions particuli\`eres doivent \^etre pr\'ealablement approch\'ees par une m\'ethode traditionnelle spectrale ou d'\'el\'ements finis par exemple. Elle est principalement d\'evelopp\'ee dans le cadre variationnel et permet la r\'esolution en temps bien plus faible que des m\'ethodes traditionnelles. Dans certain cadres, n\'eanmoins, l'approche de collocation est pr\'ef\'erable \`a l'approche variationnelle, en particulier lorsque la physique est complexe. La m\'ethode de 
 collocation r\'eduite r\'ecemment introduite  dans  \cite{ChenGottlieb} permet de poser le probl\`eme de cette fa\c con. Ainsi lorsque la m\'ethode traditionnelle est de type spectrale collocation o\`u, apr\`es avoir d\'efinit un op\'erateur discret $\mathbb{L}_\calN (\mu)$,
 on cherche un polyn\^ome $u^\calN_\mu$ tel que
$[\mathbb{L}_\calN (\mu)\, u^\calN_\mu] (x_j) = f(x_j; \mu) $ 
est v\'erifi\'e exactement sur un ensemble de   $\calN$ points de collocation $C^\calN = \{x_j\}_{j=1}^\calN$, l'approche de collocation r\'eduite propose une approximation
$u^{(N)}_{\mu^*} $:
$u^{(N)}_{\mu^*} = \sum_{j=1}^{N} c_j(\mu^*) u^\calN_{\mu^j}$ v\'erifiant (\ref{eq:rbmpde}) soit au sens des moindre carr\'e (LSRCM), puisqu'il y a plus de point $x_k$ que de coefficients $c_j$ (en effet $N<\!<\calN$) soit seulement en certain points bien choisis $ x \in C_R^N$ (ERCM).

Les m\'ethodes de collocation sont connues pour \^etre moins stables que les m\'ethodes variationnelles. Pour rectifier cet inconv\'evient, nous proposons deux types de pr\'econditionnement analytique, bas\'es sur la d\'efinition d'un op\'erateur de pr\'conditionnement $P$ et  sur une approximation de $P \mathbb{L}_\calN (\mu)\, u^\calN_\mu (x_j) \simeq P f(x_j; \mu)$ dans les deux sens pr\'ec\'edents. Les deux op\'erateurs de pr\'econditionnement analytiques sont : une version ind\'ependante du param\`etre  $P^{\mu^c} :=\mathbb{L}_\calN(\mu^c)^{-1}$, qui am\'eliore l'approximation surtout au niveau de la valeur barycentrale $\mu^c$ et une version param\'etr\'ee qui, dans le cas o\`u l'ensemble des param\`etres est le carr\'e $[0,1]^2$ repose sur une interpolation $Q_1$ entre les 4 valeurs de  $\mathbb{L}_\calN(\mu)^{-1}$ aux coins du domaine param\'etrique: 
$P^I(\mu) = P_{00} (1 - \mu^1) (1 - \mu^2) + P_{01} (1 - \mu^1) \mu^2 + P_{10} \mu^1 (1 - \mu^2) + P_{11} \mu^1 \mu^2$. Les illustrations num\'eriques des performance de ces deux pr\'econditionnement analytiques sont propos\'ees dan les figures 2 et 3. La figure 2 illustre la comparaison des trois op\'erateurs analytiques de pr\'econditionnement : sur la gauche sont trac\'es les indices d'effectivit\'e de l'estimation de l'erreur 
ÊÊ sur le syst\`eme avec ces op\'erateurs de pr\'econditionnement. Sur la droite sont trac\'es
ÊÊÊÊ la pire des convergences selon ces  sc\'enarii.  La figure 3 illustre  l'histoire de la convergence selon les op\'erateurs analytiques de pr\'econditionnement pour l'approche des moindres carr\'es (\`a gauche) et l'approche empirique de collocation (\`a droite). 

Enfin une extension de l'approche de collocation r\'eduite empirique aux cas d'EDP non lin\'eaire est aussi propos\'ee et consiste naturellementt en la v\'erification de l'EDP non lin\'eaire en des points de collocation choisis de fa\c con empirique.

\selectlanguage{english}

\vspace{-0.7cm}

\section{Introduction}
\vspace{-0.3cm}
The 
Reduced Basis Method (RBM)\cite{Barrett_Reddien,Noor_Peters,Peterson,Prudhomme_Rovas_Veroy_Maday_Patera_Turinici}
has been developed 
to numerically solve PDEs in scenarios that require a large number of numerical solutions to a parametrized
PDE, and in which we are
ready to expend significant computational time to pre-compute data that can be later used to compute accurate
solutions in real-time.
The RBM splits the solution procedure into two parts: an offline part where 
a greedy algorithm is utilized to  judiciously select $N$ parameter values for pre-computation;
and an online part where the solution for many new parameter is efficiently 
approximated by a Galerkin projection onto the low-dimensional space spanned by these $N$
pre-computed solutions.

While Galerkin methods (that are mostly used for RBM) are derived by requiring that the projection
of the residual onto a prescribed space is zero, collocation methods require the residual to be zero at some
pre-determined collocation points. 
They are attractive for their ease of
implementation, particularly for 
 time-dependent nonlinear problems \cite{HesthavenGottlieb2007,TrefethenSpecBook}.
In \cite{ChenGottlieb}, two of the authors developed 
the so-called Reduced Collocation Method (RCM). It adopts
the RBM idea for collocation methods providing a 
strategy to practitioners who prefer
a collocation, rather than Galerkin, approach. 
Our current implementation of this new method uses
collocation for both the truth solver and the online reduced solver, but the offline part could be based on a variational approach as well.
 Furthermore, one of the two approaches in \cite{ChenGottlieb}, the empirical reduced collocation
method (ERCM) allows to eliminate a potentially costly online procedure that is needed for non-affine
problems with a Galerkin approach. 
The method's efficiency 
matches  (or, for non-affine problems, exceeds) that of the traditional
RBM in the Galerkin framework.

However, collocation methods may suffer from bad conditioning.
In this paper, we propose and test two analytical  preconditioning strategies to address this issue 
in the parametric setting of RCM. 
One strategy is parameter independent, which is advantageous for ease of implementation. 
The other one is parameter dependent, but  has the traditional affinity with respect to the 
parameters which allows  extremely efficient implementation through an  offline-online decomposition.
Overall, we show that the preconditioning uniformly improves the quality of the approximation,
 and speeds up the convergence of the solution process
without adversely impacting the efficiency of the method in any significant way.
\nl{While the focus and novelty of this paper is primarily the design of the analytical  preconditioners, 
we also describe the extension of  the RCM to the nonlinear case. }
In Section \ref{sec:alg}, we briefly review 
RCM and describe our 
analytical  preconditioners. Numerical results are provided in Section \ref{sec:numerical}.

\vspace*{-0.25in}
\section{The Algorithms}
\label{sec:alg}
\vspace*{-0.15in}

We begin with  a linear parametrized PDE 
depending on a parameter $\mu \in \calD\subset \R^d$,
 written in a strong form as
$[\mathbb{L} (\mu)\, u_\mu] (x) = f(x; \mu), \,\, x \in \Omega \,\,\,\rev{\subset \mathbb{R}^n}$ 
with appropriate boundary conditions. 
We approximate the solution to this equation
using a collocation approach: for any $\mu \in \calD$, we define a discrete differentiation operator $\mathbb{L}_\calN (\mu)$, and a discrete (polynomial) solution
$u^\calN_\mu$
such that 
$[\mathbb{L}_\calN (\mu)\, u^\calN_\mu] (x_j) = f(x_j; \mu) $ 
on  a given set of  collocation points $C^\calN = \{x_j\}_{j=1}^\calN$,  usually
taken as a tensor product of $\calN_x$ collocation points for each dimension that is allowed by rectangular domains. 
We assume that  the resulting approximate solution $u^\calN_\mu$ 
is highly accurate and refer to it as the ``truth
approximation''.

\vspace{6pt}
\noindent{\em 2.1 Online algorithms}
\vspace{6pt}

For completeness, we briefly outline the RCM \cite{ChenGottlieb}. 
 The idea is that when the solution for any parameter value $\mu^* \in \calD $  is
needed, instead of solving for the costly truth approximation $u^\calN_{\mu^*}$, we 
somehow combine $N$ pre-computed truth approximations 
$u^\calN_{\mu^1},\,  \dots,\, u^\calN_{\mu^N}$  to produce a surrogate
solution $u^{(N)}_{\mu^*} $:
$u^{(N)}_{\mu^*} = \sum_{j=1}^{N} c_j(\mu^*) u^\calN_{\mu^j}.$
The key feature of the algorithm is the requirement that the surrogate solution $u^{(N)}_{\mu^*} $ 
will  satisfy the discretized differential equation in some sense
$\mathbb{L}_\calN(\mu^*) u^{(N)}_{\mu^*}  \approx f(\nl{\cdot}; \mu^*).$
Exploiting the linearity of the operator, we observe that the
system of equations we \rev{wish} to solve is : find $\vec{c}(\mu^*)$ such that 
\vspace*{-0.1in}
\begin{equation}
\vec{c}(\mu^*) = (c_1(\mu^*), c_2(\mu^*), \dots, c_N(\mu^*))^T,\quad 
\sum_{j=1}^{N} c_j(\mu^*) [\mathbb{L}_\calN(\mu^*) u^\calN_{\mu^j}](x_k)   \approx  f(x_k; \mu^*) \; \; \; \; 
.
\label{eq:rbmpde}
\end{equation}
\vspace*{-0.21in}

Satisfying the above equation exactly for $k=1,  . . ., \calN$ is usually an overdetermined system since we have  only $N$ unknowns, but $\calN >> N$ equations.

\noindent{\bf Least squares approach.}  When faced with an overdetermined system,
we  can determine the coefficients by satisfying the equation
\eqref{eq:rbmpde} in a least squares sense: 
we define, for any $\mu^*$, an $\calN \times N$ matrix
$\mathbb{A}_N (\mu^*)$ with $j^{\rm th}$ column  $\mathbb{L}_\calN(\mu^*)\, u_{\mu^j}^\calN$,
and vector of length $\calN$,
$\boldf^\calN_j (\mu^*)= f(x_j ; \mu^*) \; \; x_j \in C^\calN,$
and solve \rev{the least squares problem}
$\mathbb{A}_N^T(\mu^*)\,\mathbb{A}_N(\mu^*) \,\vec{c}(\mu^*) = \mathbb{A}_N^T(\mu^*)\, \boldf^\calN(\mu^*)$
  to obtain $\vec{c}(\mu^*)$.
 
\noindent{\bf Reduced Collocation approach.} Our second approach is more natural from the collocation point-of-view.
We determine the coefficients $\vc(\mu^*)$ by
enforcing \eqref{eq:rbmpde} at a reduced set of collocation points $C_R^N$. In other words, we solve
$\sum_{j=1}^{N} c_j(\mu^*) [ \mathbb{L}_\calN(\mu^*) u^\calN_{\mu^j}] (x)  = f(x;\mu^*)$, for $ x \in C_R^N$, 
which can also be written as $\mathbb{\rev{I}}_\calN^N \mathbb{A}_N(\mu^*) \,\vec{c}(\mu^*) = \mathbb{\rev{I}}_\calN^N  \boldf^\calN$, where $ \mathbb{\rev{I}}_\calN^N$ is a $N \times \calN$ matrix, that extracts the $N$ values of a $\calN-$vector associated to the indices of the reduced set of collocation  points.
Later we will demonstrate how this set of points can
be determined, together with the choice of basis functions, through the greedy algorithm.

\vspace{4pt}
\noindent{\em 2.2 Offline-online decomposition}
\vspace{4pt}

The size of the matrix we need for solving $\vec{c}(\mu^*)$ for each new $\mu^*$ is $N \times N$, but its assembly is 
not obviously independent of $\calN$.
For that purpose, we need the
affine assumption\footnote{$\mathbb{L} (\mu)$ can be written
as a linear combination of parameter-dependent coefficients and parameter-independent operators:
$\mathbb{L} (\mu) = \sum_{q = 1}^{Q_a} a^{\mathbb{L}}_q(\mu) \mathbb{L}_q.$
Similarly, for $f$:
$f(x; \mu) = \sum_{q = 1}^{Q_f} a^f_q(\mu) f_q (x).$} on the operator as in the Galerkin framework. Thus, the overall online component is independent of $\calN$ after a preparation 
stage where all the parameter independent quantities are precomputed 
\cite{ChenGottlieb}.  We remark that there are remedies available when the parameter-dependence is not affine \cite{Barrault_Nguyen_Maday_Patera}.

\vspace{4pt}
\noindent{\em 2.3 Analytical Preconditioning}
\vspace{4pt}

Collocation methods are frequently ill-conditioned. The situation is exacerbated when we form the normal 
equation 
in the Least Squares approach. In this section, we propose two  analytical  preconditioning techniques.  One is parameter-independent and the other is parameter-dependent.
Both will provide an operator $P$ such that 
the discretization is based on the minimization of
$P \mathbb{L}_\calN (\mu)\, u^\calN_\mu (x_j) - P f(x_j; \mu)$ 
and the reduced problem in, e.g.  the second approach   
becomes $\sum_{j=1}^{N} c_j(\mu^*) \mathbb{\rev{I}}_\calN^N \left(P \mathbb{L}_\calN(\mu^*) u^\calN_{\mu^j}\right)  = \mathbb{I}_\calN^N (P f).$

{\em Parameter-independent approach:} 
We propose using $P^{\mu^c} :=\mathbb{L}_\calN (\mu^c)^{-1}$ as a preconditioning operator. 
Here $\mu^c$ is the center of the parameter domain $\calD$.
We remark that  
this preconditioner  is in general, ideal for $\mu = \mu^c$ (making condition number exactly $1$).
Moreover, it is affordable in the parametric setting since we can perform the $\calN-$dependent operations 
for the offline preconditioning  once for all.

\vspace{-0.01in}
{\em Parameter-dependent approach:}
$P^{\mu^c}$ works well. However, it is more effective when $\mu$ is close to $\mu^c$. To have a preconditioning operator working well uniformly on the parameter domain, 
we need a parameter-dependent one. In addition, for the preconditioning to be meaningful in our parametric setting, a key requirement is that it satisfies an affine property 
similar to those for the operator $\mathbb{L}_N(\mu)$. 

Assuming our ($d-$dimensional) parameter domain is rectangular, we form $2^d$ operators at the vertices of the domain:
$ \mathbb{L}_\calN(\mu^{V_i})$ for $i = 1, \dots, 2^d$, find their inverses $P_{V_i} =  \left(\mathbb{L}_\calN(\mu^{V_i}) \right)^{-1}$, and  define the preconditioning operator through interpolation. In the case $d=2$ (we assume $\mu = (\mu^1,\mu^2) \in [0,1]^2$ without loss of generality), this is a $Q_1$ interpolation defined as below:
$P^I(\mu) = P_{00} (1 - \mu^1) (1 - \mu^2) + P_{01} (1 - \mu^1) \mu^2 + P_{10} \mu^1 (1 - \mu^2) + P_{11} \mu^1 \mu^2,$
where $P_{ij}$ is $P_{V} $ with $V$ being the $(i,j)-$corner.

\vspace{4pt}
\noindent{\em 2.4 Offline algorithms}
\vspace{4pt}

In this section we describe 
 the two greedy algorithms for the least squares and the reduced collocation approaches
 for choosing the reduced basis set $\{ u_{\mu^1}^\calN, \dots,
u_{\mu^N}^\calN \}$. We assume that  given   $\{ u_{\mu^1}^\calN, \dots,
u_{\mu^i}^\calN \}$ we can compute  an upper bound $\Delta_i(\mu)$  for the error of the reduced solution
$u_\mu^{(i)}$ for any parameter $\mu$ \cite{ChenGottlieb}.

\begin{algorithm}[htp]
  \caption{Least Squares Reduced Collocation Method (LSRCM)\rev{: Offline Procedure}}\label{alg:LSgreedy}
  \begin{algorithmic}
\STATE {\bf 1.} Discretize the parameter domain $\calD$ by $\Xi$, and denote the center of $\calD$ by $\mu^c$.
\medskip
\STATE {\bf 2.} Randomly select $\mu^1$, solve $\mathbb{L}_\calN (\mu^1)\, u^\calN_{\mu^1}
(x) = f(x; \mu^1)$ for $x \in C^\calN$ and let $\xi_1^\calN = u^\calN_{\mu^1}$.
\medskip
\STATE {\bf 3.} For $i = 2, \dots, N$ do
\begin{itemize}
\item [{\bf i).}] For all $\mu \in \Xi$, form $\mathbb{A}_{i-1}(\mu) = \left(\mathbb{L}_\calN(\mu)\, \xi_{1}^\calN, \mathbb{L}_\calN(\mu)\, \xi_{2}^\calN, \, \dots, \,
\mathbb{L}_\calN(\mu)\, \xi_{{i-1}}^\calN\right)$.
\item [{\bf ii).}] For all $\mu \in \Xi$, solve
$\mathbb{A}_{i-1}(\mu)^T\,\mathbb{A}_{i-1}(\mu) \,\vec{c} = \mathbb{A}_{i-1}^T(\mu)\,\boldf^\calN$ to obtain $u^{(i-1)}_\mu
= \sum_{j = 1}^{i - 1} c_j \xi^\calN_{j}$ and $\Delta_{i-1}(\mu)$.
\item [{\bf iii).}] Set $\mu^i = argmax_{\mu}\,\,\Delta_{i-1}(\mu)$, and solve $\mathbb{L}_\calN (\mu^i)\, u^\calN_{\mu^i}(x) = f(x; \mu^i)$ for $x \in C^\calN$.
\item [{\bf iv).}] Apply a modified Gram-Schmidt transformation, with inner product defined by \\$(u,v) \equiv \left(\mathbb{L}_\calN (\mu^c)u, \mathbb{L}_\calN (\mu^c) v\right)_{L^2(\Omega)}$, on $\left\{\xi^\calN_{1}, \xi^\calN_{2}, \dots, \xi^\calN_{i-1}, u^\calN_{\mu^i}\right\}$  to obtain 
$\left\{\xi_1^\calN, \xi_2^\calN, \dots, \xi_i^\calN\right\}$.
\end{itemize}
  \end{algorithmic}
\end{algorithm}

\vspace*{-0.1in}

\begin{algorithm}[h!]
  \caption{Empirical Reduced Collocation Method (ERCM)\rev{: Offline Procedure}}\label{alg:RCgreedy}
  \begin{algorithmic}
\STATE {\bf 1.} Randomly select $\mu^1$, solve $\mathbb{L}_\calN (\mu^1)\, u^\calN_{\mu^1}
(x) = f(x; \mu^1)$, let $x^1 = argmax_{x \in X} \,\, \left|u^\calN_{\mu^1}(x)\right|,\quad \xi^\calN_1 =
\frac{u^\calN_{\mu^1}}{u^\calN_{\mu^1}(x^1)}$.
\STATE {\bf 2.} For $i = 2, \dots, N$ do
\begin{itemize}
\item [{\bf i).}] Let $C_R^{i-1} = \left\{x^1,\dots,x^{i-1}\right\}$.
\item [{\bf ii).}] For all $\mu \in \Xi$, solve \ $\sum_{j=1}^{i-1} c_j \mathbb{I}_\calN^N
\left(\mathbb{L}_\calN(\mu) u^\calN_{\mu^j}\right)  = f(x; \mu) \,\, {\rm for} \,\, x \in C_R^{i-1}$ to
obtain $u^{(i-1)}_\mu = \sum_{j = 1}^{i - 1} c_j u^\calN_{\mu^j}$.
\item [{\bf iii).}] Set $\mu^i = argmax_{\mu \in \Xi}\,\,\Delta_{i-1}(\mu)$ and solve $\mathbb{L}_\calN (\mu^i)\, u^\calN_{\mu^i} (x) = f(x; \mu^i)$.
\item [{\bf iv).}] Find $\alpha_1, \dots, \alpha_{i-1}$ such that, if we define $\xi^\calN_i = u^\calN_{\mu^i} - \sum_{j=1}^{i-1}
\alpha_j\,\xi^\calN_j$, we have $\xi^\calN_i(x^j) = 0$ for $j = 1, \dots, i-1$.
\item [{\bf v).}] Set $x^i = argmax_x \,\, \left|\xi^\calN_{i}\right|$ and $\xi^\calN_i = \frac{\xi^\calN_i}{\xi^\calN_i(x^i)}$.
\item[\rev{\bf vi).}] \rev{Apply modified Gram-Schmidt transformation on $\left\{\xi_1^\calN, \dots, \xi_i^\calN\right\}$.}
\end{itemize}
\medskip
  \end{algorithmic}
\end{algorithm}

\vspace{3pt}
\noindent{\em 2.5 Extension to the nonlinear case}
\vspace{3pt}

\nl{The ERCM approach is more economical than the EIM implementation of the 
variational RBM for linear problems having a large number of varying coefficients, such as 
the case when geometry is a parameter  \cite{lovgren2011reduced}.
Here, we outline the procedure when we have a general nonlinear operator 
$\mathbb{G}(u; \mu) + \mathbb{L}(\mu) u$ where  $\mathbb{G}(u; \mu)$ is nonlinear. 
The parameter dependence can be handled in the same way with possibly the Empirical Interpolation \cite{Barrault_Nguyen_Maday_Patera} needed, so it suffices to assume $\mathbb{G}(u; \mu) \equiv \mathbb{G}(u)$. 
In the following, we present our approach through the example of  the viscous Burgers' equation 
$u u_x - \mu u_{xx} =f_\mu $; The formal extension to multi-dimension vector equations   is straightforward. 
In the case of Burgers' equation, the discretized system for any parameter  $\mu^j$ becomes 
$\bu \odot \mD \bu - \mu^j \mD_2 \bu =f_{\mu^j},$
where $\bu$ is the $\calN \times 1$ vector containing the point values of $u$ on the Chebyshev grid, 
$\odot$ is the well-known Hadamard product for vectors that denotes element wise multiplication,
and $\mD$ and $\mD_2$ are the $\calN \times \calN$  first and second order differentiation matrices, respectively.
We assume that we are given $N$ solutions $u_{\mu^j}^\calN$ and we define $u^{(N)}_{\mu^*} = \sum_{j=1}^N c_j(\mu^*) u_{\mu^j}^\calN$.
We find the values of the unknown coefficients $c_j$ by satisfying a nonlinear equation of the form
$ G(\vc) - \mu^* \mL \vc = f_{\mu^*}$. Here,  $\vc$ is the column vector of coefficients $c_j$ of length $N$, 
$\mL$ is  an $N \times N$ matrix, and $G: \mathbb{R}^{N \times 1} \rightarrow \mathbb{R}^{N \times 1}$ 
is a nonlinear function.  These solutions are then solved by some iterative fixed-point like method. 
For the ERCM case, $\mL$ and $G$ come from the discrete solution satisfying  $\bu \odot \mD \bu - \mu^* \mD_2 \bu =f_{\mu^*}$ 
at a reduced set of collocation points $C_R^N$. 
}

\vspace{-0.25in}

\section{Numerical Results}
\label{sec:numerical}
\vspace{-0.1in}

\begin{figure}
\centering
\begin{minipage}{.32\textwidth}
  \centering
        \includegraphics[width=\textwidth]{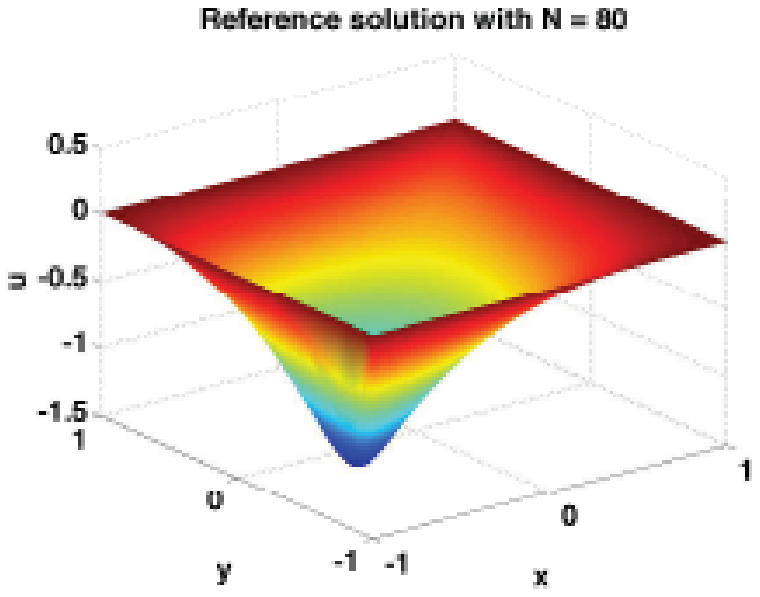}
  \caption{The truth approximation
    for $\mu_1 = 1$ and $\mu_2 = 0.5$ computed on a $81 \times 81$ Chebyshev grid.}
  \label{fig:truthsample}
\end{minipage}%
\hspace*{.04\textwidth}
\begin{minipage}{.64\textwidth}
  \centering
          \includegraphics[width=0.48\textwidth]{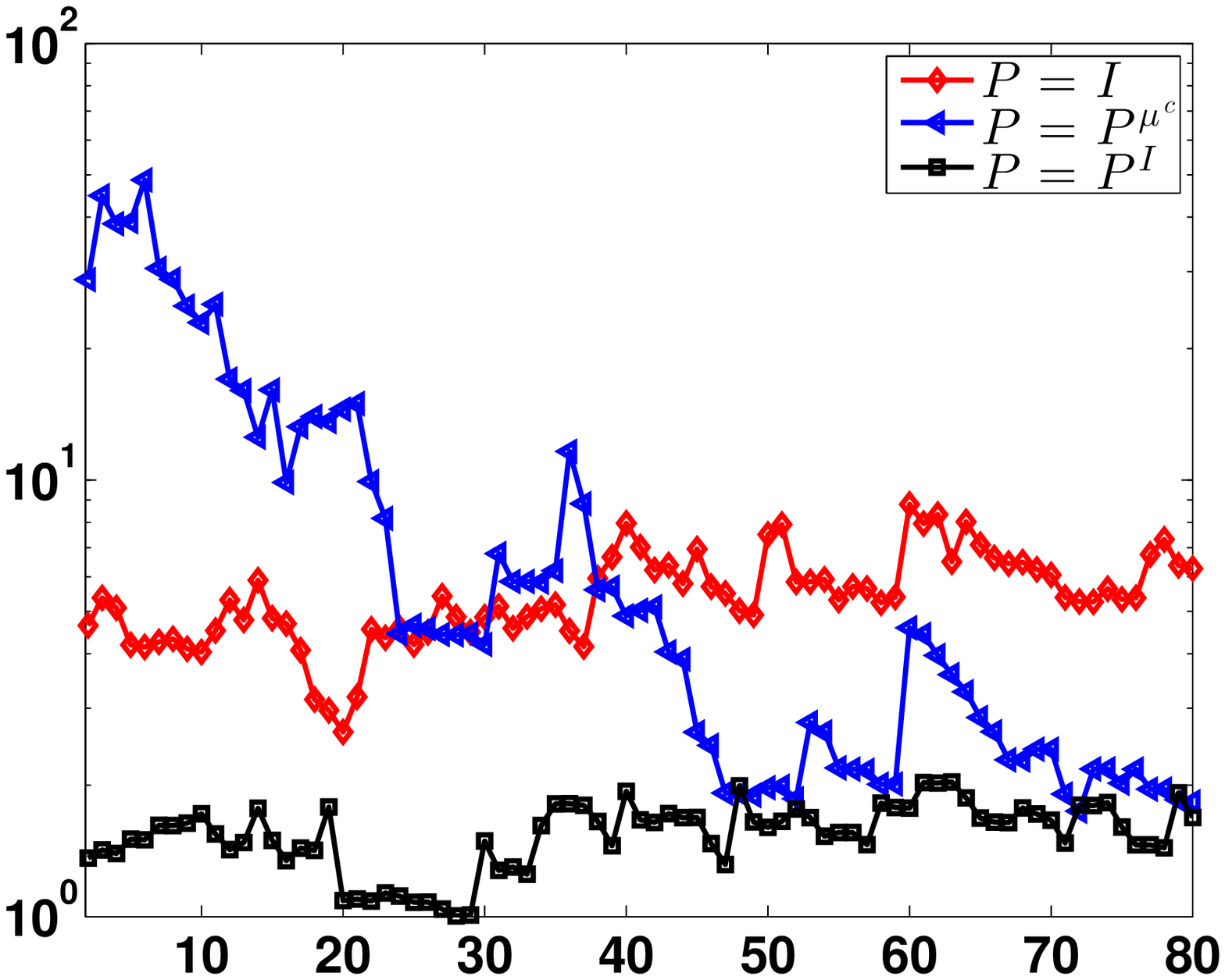}
          \includegraphics[width=0.48\textwidth]{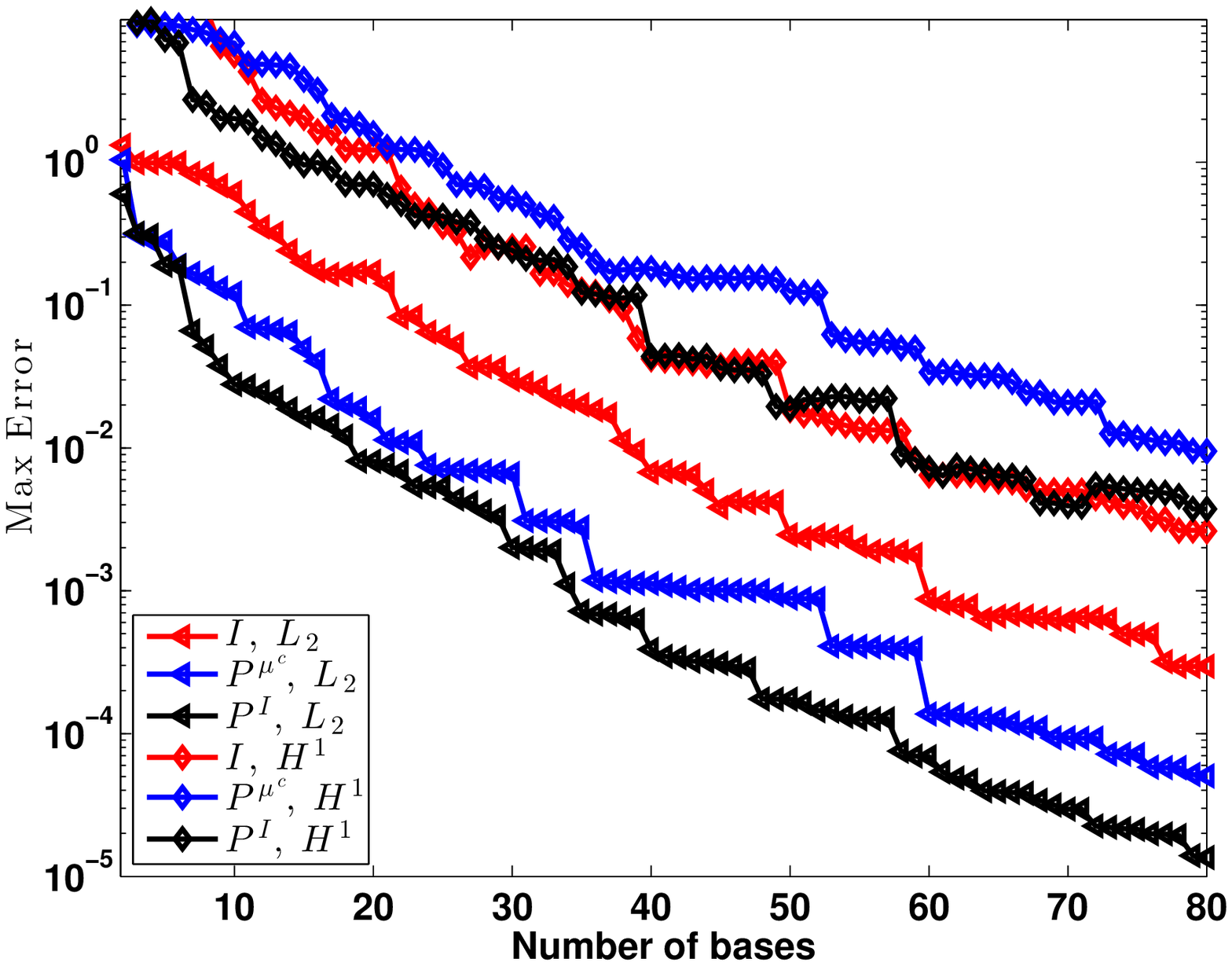}
  \caption{Comparison of the three analytical  preconditioning operators: On the left is the plot for the effectivity indices for the error estimate 
  on the system with these preconditioning operators. On the right are 
    the worst case scenario convergence plots.}
  \label{fig:L2d_H1u}
\end{minipage}
\end{figure}
\begin{SCfigure}
       \includegraphics[width=0.24\textwidth]{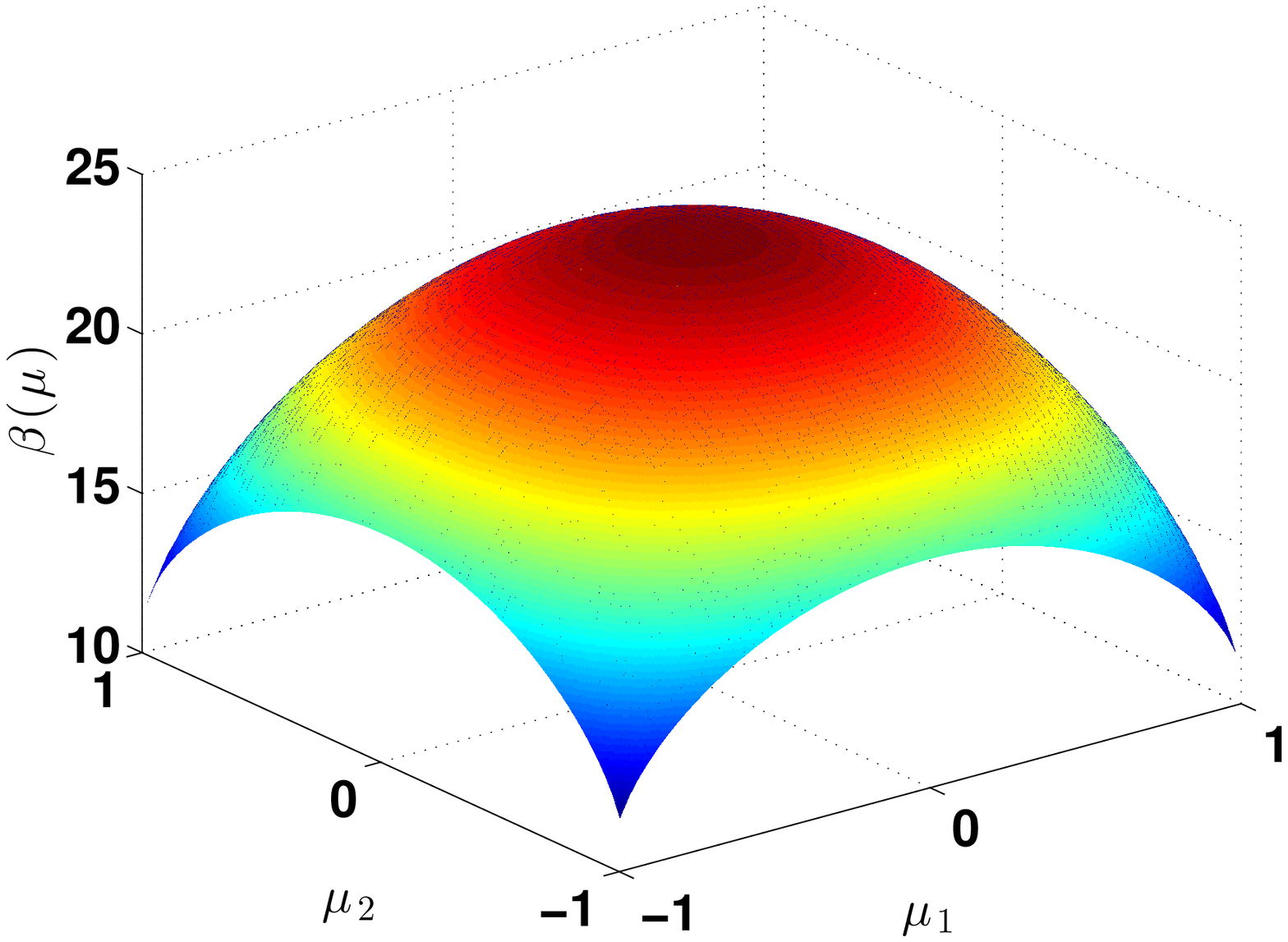}
        \includegraphics[width=0.24\textwidth]{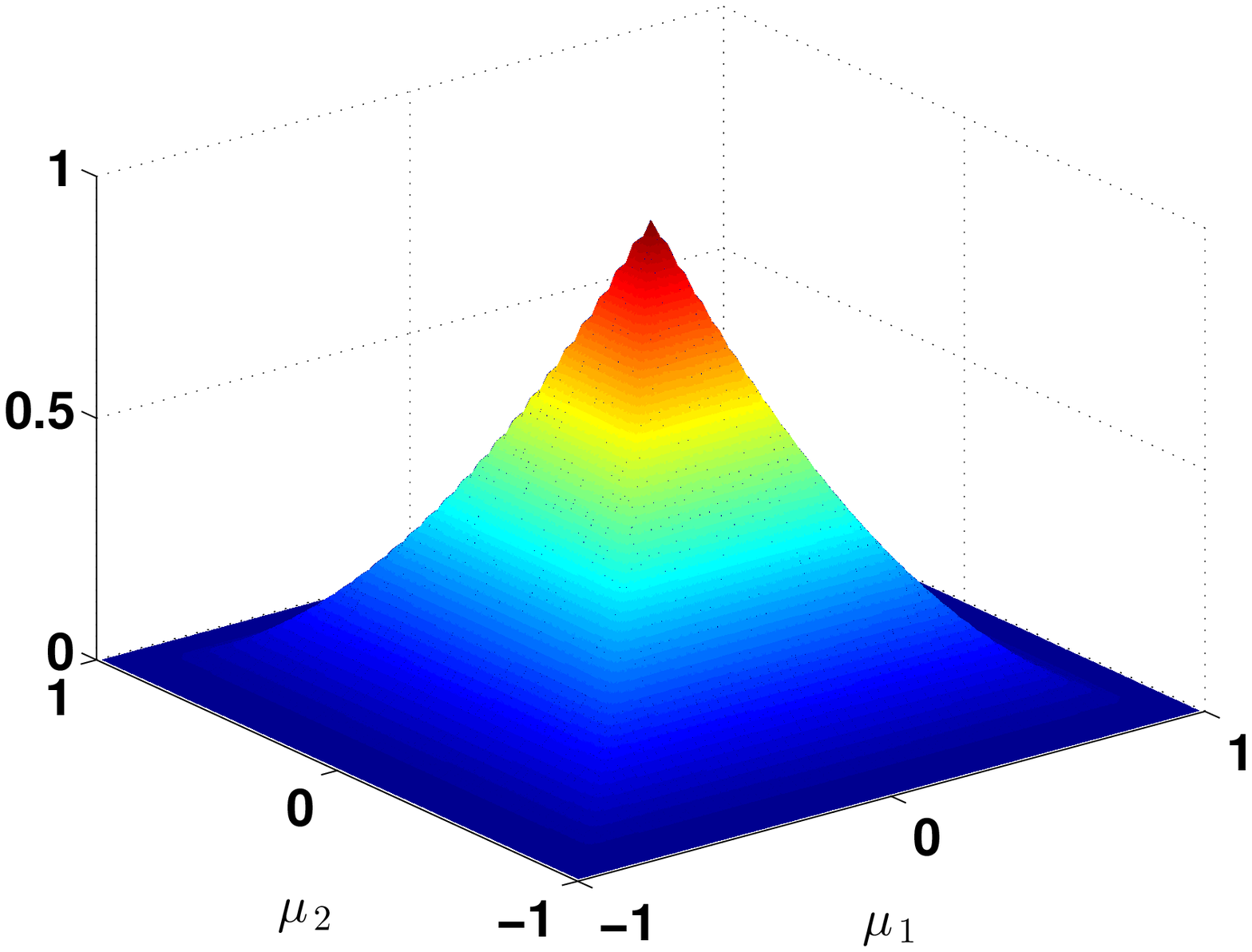}
        \includegraphics[width=0.24\textwidth]{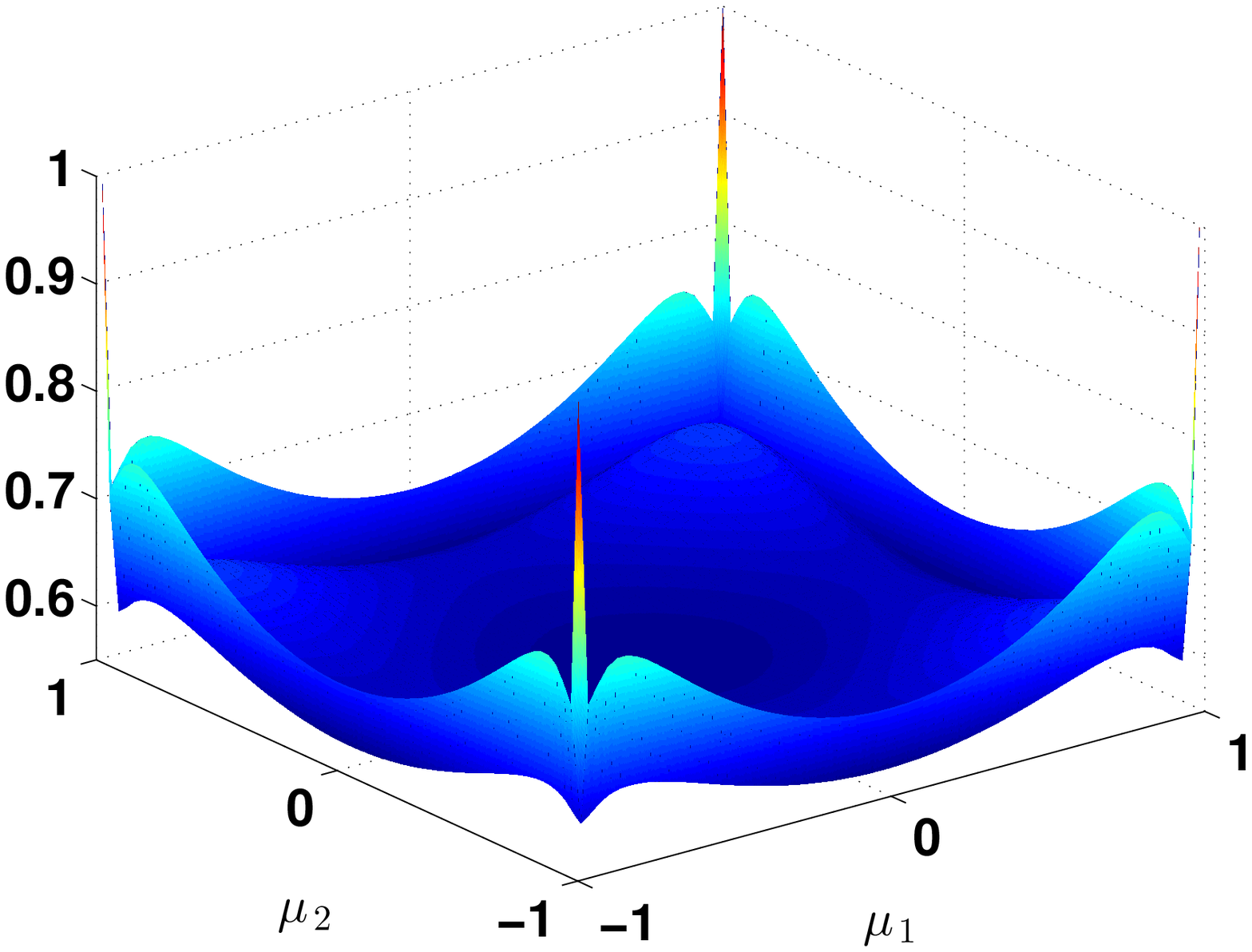}
    \caption{\small From left to right: the non-preconditioned inf-sup number ($P = I$), Parameter-independent preconditioning ($P = P^{\mu^c}$), 
    and $Q_1$-interpolating parametric preconditioning ($P = P^I(\mu)$).}
    \label{fig:howtoprec}
\end{SCfigure}
\begin{SCfigure}
        \includegraphics[width=0.35\textwidth]{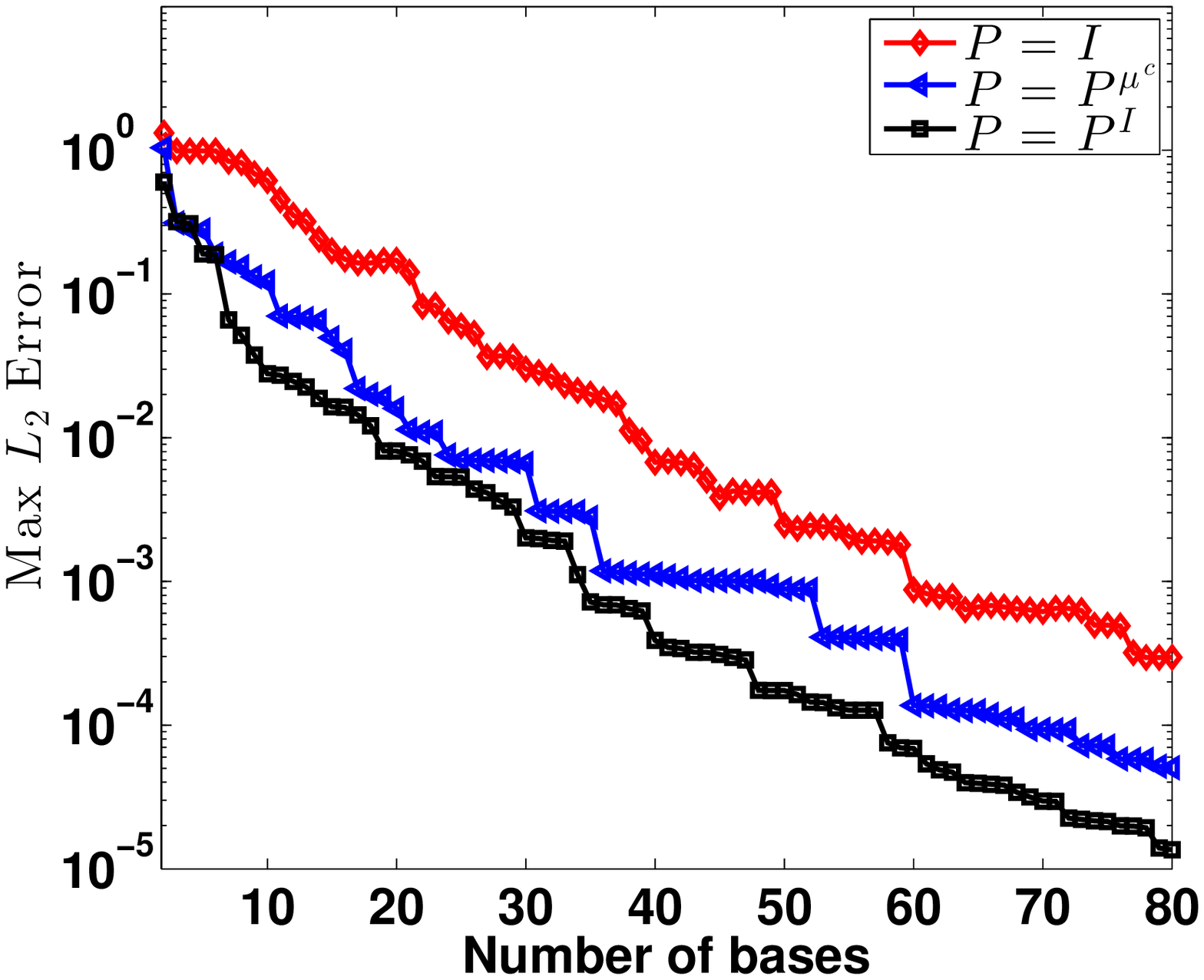}
        \includegraphics[width=0.35\textwidth]{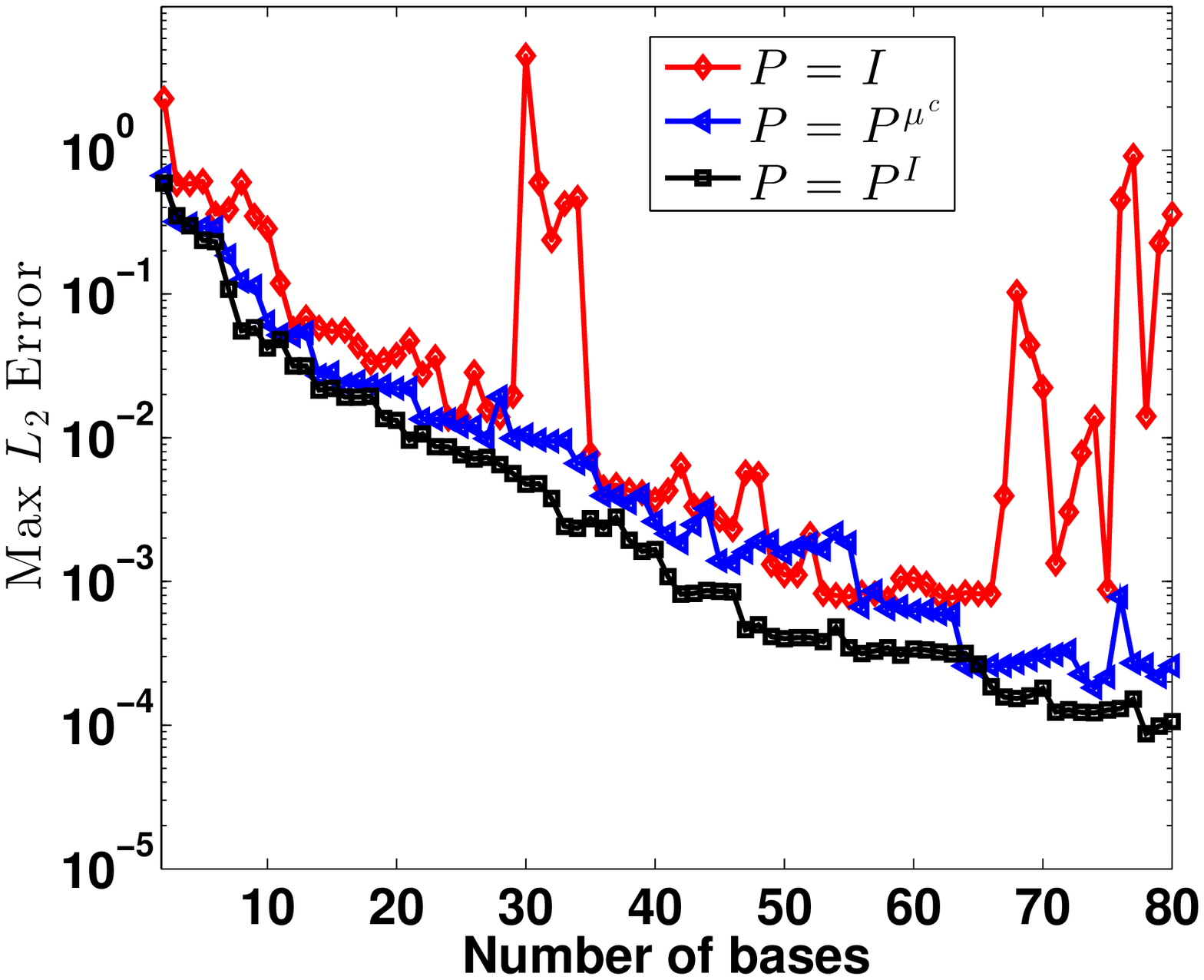}
\caption{\small Comparison of history of convergence when different analytical  preconditioning operators is used for the least squares approach (left) and empirical collocation approach (right). Maximum is taken among $L_2$ errors between truth approximation and reduced basis solution for $1,057$ randomly selected parameter values.}\label{fig:IvsFD_LS_Ellip}
\end{SCfigure}

In this section, we demonstrate  the accuracy and efficiency of
the proposed methods on a 2D  diffusion-type problem with zero Dirichlet boundary condition \cite{ChenGottlieb}:
$(1+ \mu_1 x) u_{xx} + (1+ \mu_2 y) u_{yy}= e^{4 x y} 
\,\, \mbox{ on } \Omega = [-1,1] \times [-1,1].$
Our truth approximations are generated by a spectral Chebyshev collocation method
\cite{TrefethenSpecBook,HesthavenGottlieb2007}. For $C^\calN$, we use the Chebyshev  grid based
on \rev{$\calN_x$} points in each direction \rev{with $\calN_x^2 = \calN$}. We consider the parameter domain $\calD$ for $(\mu_1,\mu_2)$  to be
$[-0.99,0.99]^2$
. For $\Xi$, they are
discretized uniformly by a $64 \times 64$ 
Cartesian grid. 
One sample solution for this problem is plotted in Figure \ref{fig:truthsample}.

In Figure \ref{fig:L2d_H1u} (Left), we show that while the non-parametric preconditioning $P^{\mu^c}$ give non-uniform improvement,  the parametric preconditioning $P^I(\mu)$ improves effectivity indices by one order of magnitude.
Figure \ref{fig:L2d_H1u} (Right) shows that the $P^{\mu^c}$ preconditioning operator
  improves the $L_2$ norm of the error but worsens the $H^1$ norm.
In comparison, $P^I(\mu)$ improves the error in $L^2$ norm without significantly degrading (in some cases improving) 
the error in $H^1$ norm. 
Finally, we plot the stability constant of these preconditioned operators as a function of the parameter 
in Figure \ref{fig:howtoprec}. 
We clearly see that $P^I(\mu)$ is most efficient in terms of enforcing the parametric stability number uniformly close to $1$. 
\nl{We also tested diagonal preconditioning $P^I_D(\mu)$ (not reported here) by using the same interpolating procedure as $P^I(\mu)$ and replacing the inverses of the full operators  by the inverses of the diagonals. 
Clearly $P^I_D(\mu)$ is cheaper to compute than the other preconditioners, but its performance is significantly worse than $P^I(\mu)$ and even worse than $P^{\mu^c}$.}

For the preconditioned RCM, we can see, from Figures \ref{fig:IvsFD_LS_Ellip} that the error for the least squares approach is around one order of magnitude better in the worst case scenario.
For the empirical reduced collocation approach, the error is smaller and, more importantly,
 converges much more stably.

\vspace{-.3in}
\section{Concluding Remarks}
\label{sec:conclude}
\vspace{-.2in}
We propose and test two analytical preconditioning strategies in the context of reduced collocation method.
The parameter dependent one is shown to be capable of offline-online decomposition, improving both
the quality of error estimation uniformly on the parameter domain, and enabling the preconditioned 
reduced collocation method to converge much faster and more stably than the non-preconditioned version.


\end{document}